\documentclass[a4paper,10pt]{amsart}
\usepackage{amsmath}
\usepackage{amssymb}
\usepackage{amsfonts}
\usepackage{amscd}
\usepackage{latexsym,graphicx,caption}
\input xy
\xyoption{all}
\usepackage{verbatim}

\newcommand{\ti}[1]{\tilde{#1}}
\newcommand{\ov}[1]{\overline{#1}}

\newcommand{\lo}{\longrightarrow}

\newcommand{\diam}{{\hfill \nobreak} $\Box$}
\newcommand{\noi}{\noindent}
\newcommand{\lspace}{\vspace{0.45cm}}
\newcommand{\sspace}{\vspace{0.1cm}}

\newtheorem{theor}{\noi Theorem}
\newtheorem*{theore}{\noi Theorem}
\newtheorem{lemma}{\noi Lemma}

\newtheorem{corol}{\noi Corollary}

\newtheorem{prop}{\noi Proposition}
\newtheorem{defi}{\noi Definition}

\newtheorem{conj}{\noi Conjecture}

\theoremstyle{remark} 
\newtheorem{rem}{\noi Remark}

\newcommand{\ZZ}{{\mathbb Z}}
\newcommand{\QQ}{{\mathbb Q}} 
 
\newcommand{\RR}{{\mathbb R}}
\newcommand{\CC}{{\mathbb C}}
\DeclareMathOperator{\Hom}{Hom}
\DeclareMathOperator{\Aut}{Aut}

\newcommand{\G}{\mathbf G}

\newcommand{\gc}{{\mathfrak g}_{\CC}}

\newcommand{\kc}{{\mathfrak k}_{\CC}}
\newcommand{\pc}{{\mathfrak p}_{\CC}}

\newcommand{\proj}{\mathbf P}

\newcommand{\Ga}{\Gamma}

\newcommand{\U}{\mathbf U}

\newcommand{\ggc}{G_{ \mathbb C}}

\newcommand{\hlc}{{\mathfrak h}_{ \mathbb C}}
\newcommand{\lieg}{{\mathfrak g}}
\newcommand{\liek}{{\mathfrak k}}
\newcommand{\lieh}{{\mathfrak h}}
\newcommand{\liet}{{\mathfrak t}}
\newcommand{\lieu}{{\mathfrak u}}
\newcommand{\liel}{{\mathfrak l}}

\newcommand{\liep}{{\mathfrak p}}
\newcommand{\liev}{{\mathfrak v}}
\newcommand{\kg}{{\mathfrak k}}

\newcommand{\vc}{{\mathfrak v}_{ \mathbb C}}
\newcommand{\n}{{\mathfrak n}}
\newcommand{\nplus}{{\mathfrak n_{+}}}
\newcommand{\nmoins}{{\mathfrak n_{-}}}

\newcommand{\q}{{\mathfrak q}}

\newcommand{\bul}{\bullet}

\newcommand{\PU}{\mathbf{PU}}
\newcommand{\PGL}{\mathbf{PGL}}
\newcommand{\GL}{\mathbf{GL}}

\newcommand{\dtheta}{\check{D}_{\theta}}

\newcommand{\HHom}{\textnormal{Hom}}
\newcommand{\Qtheta}{\mathbf{Q}_{\theta}}

\newcommand{\End}{\textnormal{End}}
\newcommand{\Ad}{\textnormal{Ad}}

\newcommand{\FF}{\mathcal{F}}
\newcommand{\rk}{{\text{rk}}}
\newcommand{\GC}{\mathbf{G}_\CC}
\newcommand{\MD}{\mathbf{M}_{\textnormal{Dol}}}

\newcommand{\holP}{\textnormal{P}}

\newcommand{\PP}{\mathcal{P}}
\newcommand{\delbar}{\overline{\partial}}
\newcommand{\cvhs}{$\CC$VHS}
\newcommand{\ad}{\textnormal{ad}}
\newcommand{\TT}{{\mathbf T}}
\newcommand{\M}{\mathbf{M}}

\newcommand{\qc}{\mathfrak{q}}
\newcommand{\Gr}{\textnormal{Gr}}
\newcommand{\Res}{\textnormal{Res}_{\CC/\RR}}
\newcommand{\ea}{e_{\alpha}}
\newcommand{\eb}{e_{\beta}}
\newcommand{\eg}{e_{\gamma}}
\newcommand{\ee}{e_{\varepsilon}}
\newcommand{\Tr}{\text{Tr }}
\newcommand{\E}{\mathbb{E}}

\title{On the cohomology of K\"ahler groups}

\author{Bruno KLINGLER}

\begin{document}
\baselineskip 14pt
\setcounter{tocdepth}{1}
\maketitle
\tableofcontents
\section{Introduction and results}

\subsection{The main result}
This is the second of two papers (c.f. \cite{kling}) investigating which
finitely generated groups can be realized as fundamental groups
of a connected compact K\"ahler manifold. Such groups are called K\"ahler
groups. Even if this question is completely open there are many known
obstructions for an infinite finitely presented group being K\"ahler (we refer to \cite{abckt}
for a panorama). Most of them come from Hodge theory (Abelian or not)
{\em in cohomological degree one}. The most interesting conjecture concerning infinite
K\"ahler groups, due to Carlson and Toledo and publicized by Kollar \cite{ko} and
Simpson, deals with {\em cohomology in degree $2$}~:
\begin{conj}[Carlson-Toledo] \label{conj1}
Let $\Ga$ be an infinite K{\"a}hler group. Then virtually $b^2(\Ga) >0$.
\end{conj}

\begin{rem} Recall that a group $\Ga$ has virtually some property $\mathcal{P}$ if
a finite index subgroup $\Ga' \subset \Ga$ has $\mathcal{P}$.
\end{rem}

Our first result in this paper is the following~:

\begin{theor} \label{main}
Let $\Ga$ be a finitely presented group admitting a morphism $\rho:
\Ga \lo GL(n, \CC)$ with unbounded image in $GL(n, \CC)$. 

If $\Ga$ is a K\"{a}hler group
then $b^2(\Ga) + b^4(\Ga) >0$.

Moreover if $b^2(\Ga)=0$ then the bounded cohomology group
$H^4_b(\Ga, \RR)$ is non-trivial.
\end{theor}

Let us make some comments on theorem~\ref{main}. 
\begin{itemize}
\item In \cite{to} Toledo constructed infinite K\"ahler groups (in fact~:
  fundamental groups of complex smooth projective varieties) which are
  not residually finite. In particular such groups are not
  linear. However all known infinite K\"ahler groups satisfy the
  assumption of theorem~\ref{main}~: they admit a morphism to some
  complex linear group $GL(n, \CC)$ with unbounded image.
\item Restricting ourselves to such groups rather than infinite
  groups is a classical drawback in the study of K\"ahler groups~: we know
  nothing about their non-linear part.
\item The main content of theorem~\ref{main} is of course that $b^2(\Ga) +
  b^4(\Ga) >0$. Although the remark concerning bounded cohomology
  might be of interest I don't know many groups $\Ga$ for which one
  can compute $H^4_b(\Ga, \RR)$.
\end{itemize}

\subsection{A variant of \cite[theor.1]{kling}}
In \cite{kling} we proved the following~:
\begin{theore} \cite[theor.1]{kling} \label{linear}
Let $\Ga$ be an infinite linear $r$-dimensional duality group, $r \geq
6$. If $\Ga$ is a K\"ahler group then virtually
$$b^2(\Ga) + b^4(\Ga) >0 \;\;.$$
\end{theore}

Theorem~\ref{main} may look more general than
\cite[theor.1]{kling}~: there is no duality nor linearity
assumption, and one does not have to pass to a finite index subgroup. However it is crucial to realize that theorem~\ref{main} 
does not imply \cite[theor.1]{kling} as there are many linear groups
whose complex representations are all bounded. Typically
\cite[theor.1]{kling} implies the following~:
\begin{theore} \cite[theor.2]{kling} \label{p-adic}
Let $G_v$ be the group of $F_v$-points of an algebraic
group $\G_v$, with reductive neutral component, over a non-Archimedean
local field $F_v$ of characteristic $0$. 
Suppose that $\text{rank}_{F_{v}} \G_v \geq 6$.

Then a cocompact lattice $\Ga \subset G_v$ is never the
fundamental group of a smooth projective complex variety.
\end{theore}

On the other hands theorem~\ref{main} says nothing in this situation~:
such a $\Ga \subset G_v$ is linear but every representation $\rho :
\Ga \lo GL(n, \CC)$ has bounded image by Margulis's superrigidity
theorem \cite{mar}.

\sspace
It is better to think of theorem~\ref{main}, respectively \cite[theor.1]{kling}
as the geometric facet, respectively the topological facet, of the same phenomenon. As we will see the main ingredient for
theorem~\ref{main} is non-Abelian Hodge theory and a careful
differential-geometry study of Pontryagin forms on period domains.  On
the other hand \cite[theor.1]{kling} relies on the topological
properties of Stein spaces and duality properties in group cohomology. 

\sspace
In fact putting together theorem~\ref{main} and the proof of
\cite[theor.1]{kling} we obtain the following improvement of
theorem~\ref{main} for linear duality groups~:

\begin{theor} \label{linear2}
Let $\Ga \subset GL(n, \CC)$ be an infinite linear $r$-dimensional duality group, $r \geq
6$, with $b^2(\Ga)=0$. Suppose that $\Ga$ is not bounded in $GL(n,
\CC)$. If $\Ga$ is K\"ahler then virutally $b^4(\Ga)>1$.
\end{theor}

\subsection{Application to lattices in real Lie groups}

A simple way of constructing ``big'' K\"ahler groups (and essentially
the only known one) consists in generalizing to higher dimensions the
uniformization of a smooth projective complex curve of genus $g \geq
2$ as quotient of the Poincar\'e disk. Let $G$ be a real simple Lie
group of Hermitian non-compact type, i.e. such that the associated
symmetric space $G/K$ admits an unbounded $G$-invariant K\"ahler
metric. For example $G = PU(n,1)$ the group of isometries of the unit
ball $\mathbf{B}^n_\CC = \{[z_1, \cdots z_{n+1}] \in \proj^n\CC, \;\;
|z_1|^2 + \cdots + |z_n|^2 < |z_{n+1}|^2 \}$. Any torsion free
cocompact lattice $\Ga$ of $G$ acts properly discontinuously on $G/K$ and the
quotient $X = \Gamma \backslash G/K$ is a compact K\"ahler manifold
(in fact a smooth projective complex variety) with fundamental group
$\Ga$. One can construct more sophisticated examples starting from these ones
using Lefschetz's hyperplane theorem and glueing techniques \cite{to}.

In this context theorem~\ref{linear2} enable us to prove (with some
restrictions on the rank) another conjecture by
Carlson-Toledo (certainly the main motivation for conjecture~\ref{conj1})~:
\begin{theor} \label{lattice}
Let $G$ be a simple real Lie group of non-compact type, of real rank
at least $20$, and $\Ga$ a cocompact lattice in $G$. Then $\Ga$ if a
K\"ahler group if and only if $G$ is of Hermitian type.
\end{theor}

\begin{rem}
In \cite{si1} Simpson showed that if a cocompact lattice $\Ga$ in a
simple real Lie group $G$ is K\"ahler then $G$ is necessarily of
Hodge type (i.e. admits discrete series). Thus theorem~\ref{lattice} consists in reducing the class
of groups of Hodge type to the class of groups of Hermitian type. That is~:
excluding the two families $G= SO(2p, q)$, $p \not =1$ and $G = Sp(p,q)$,
$p \geq q$, $q \geq 1$ and eight exceptional groups. Our proof does
not apply in small rank. Certainly the most interesting remaining case
is the following~: does there exists a cocompact lattice in $Sp(n,1)$
which is a K\"ahler group~?
\end{rem}

\subsection{Ingredients for the proof}
Let $\Ga$ be the fundamental group of a compact connected Riemannian manifold $X$ and $\rho: \Ga \lo G=\G(\RR)$ be a
Zariski-dense representation, where $\G$ denotes a real almost simple
algebraic group (non-anisotropic). One can try (c.f. \cite{kling2}) to understand $\rho$ by
studying the
pull-back map in cohomology
$$\rho^* :H^\bullet_{ct}(G, \RR) \lo H^\bullet(\Ga, \RR)\;\;.$$

The continuous cohomology $H^\bullet_{ct}(G, \RR)$ is completely
understood since the $60's$. It is naturally isomorphic to the De Rham
cohomology $H^\bullet(A^\bullet(G/K)^G,d)$ of $G$-invariant
differential forms on the symmetric space $G/K$ or equivalently with the
relative Lie algebra cohomology $H^\bullet(\lieg, K, \RR)$.
On the other hand the datum of $(X, \rho)$ essentially uniquely
defines an harmonic $\rho$-equivariant map
$$\overline{f_\rho} : \ti{X} \lo G/K\;\;,$$ 
where $\ti{X}$ denotes the universal cover of $X$.
Contemplating the commutative diagram
\begin{equation} \label{diag}
\xymatrix{
H^\bullet_{ct}(G, \RR) \ar@{=}[d] \ar[r]^{\rho^{*}} & H^\bullet(\Ga, \RR)  \ar[d]  \\
H^\bul(A^*(G/K)^G, d) \ar[r]_>>>>>>{\overline{f_\rho}^* } & H^\bullet(X, \RR)
}
\end{equation}
we would like to exhibit some class $c \in H^\bullet_{ct}(G, \RR)$ whose
image $\overline{f_{\rho}}^*(c)$ does not vanish in $H^\bullet(X,
\RR)$, thus showing that the class $\rho^*(c) \in H^\bullet(\Ga, \RR)$ does not
vanish.
For a general $X$ not much can be done. When $X$ is a compact K\"ahler
manifold this picture can be much improved using non-Abelian Hodge
theory.

\sspace
Let $X$ and $\rho$ be as in theorem~\ref{main}. If $\Ga$ admits a
non-rigid complex finite dimensional representation then Simpson's
version of Lefschetz's
theorem readily implies that $b^2(\Ga)>0$
(c.f. \cite[Appendix B]{kling}). Thus
we can assume that $\rho$ is reductive rigid. As $\rho$ is unbounded non-Abelian
Hodge theory implies that $\rho$ is the monodromy of a (polarized) non-trivial
complex variation of Hodge structure (\cvhs\- for short) on $X$~: there
exists a $C^\infty$ complex vector bundle $E$ (say of rank $n$) with a
decomposition $E=\oplus_{i=0}^k E^i$ (called the Hodge decomposition) into a direct sum of $C^\infty$
subbundles, a flat connection $\nabla$ with monodromy $\rho$
satisfying Griffiths's transversality condition 
$$\nabla:E^i\lo A_X^{0,1}(E^{i+1})\oplus A_X^{1,0}(E^i)\oplus A_X^{0,1}(E^i)\oplus A_X^{1,0}(E^{i-1})
$$
and a parallel Hermitian form $h$ which makes the Hodge decomposition
orthogonal and is positive on $E^i$ if $i$ is even and negative if $i$
is odd. If we set $r_i=\rk E^i$, the monodromy representation $\rho$
has image contained in a group isomorphic to
$U(p,q)\subset GL(n=p+q,\CC)$ where $p=\sum_{i\ {\rm even}}r_i$ and
$q=\sum_{i\ {\rm odd}}r_i$. 

The Hodge filtration of $E$ by the subbundles $F^s=\oplus_{i\geq s}
E^i$ is holomorphic. Let us fix $x_0\in X$ and let
$\E=\oplus_{i=0}^k\E^i$ be the fiber of $E$ over $x_0$ endowed with
the restriction $h_0$ of the Hermitian metric $h$ on $E$. Let $\GC$ be
the complex reductive algebraic group Zariski-closure of $\rho(\Ga)$ in $\PGL(\E)$, let $\G$ be its real form
$\GC \cap\PU(\E,h_0)$, let $Q$ be the parabolic
subgroup of $G_\CC=\G_\CC(\CC)$ which stabilizes the flag ${\mathbb F}^k\subset
{\mathbb F}^{k-1}\subset\dots\subset {\mathbb F}^1\subset {\mathbb
  F}^0=\E$ defined by the Hodge filtration and let $V=G\cap Q\subset{\rm
  P}(U(\E^0,{h_0}_{|\E^0})\times\dots\times U(\E^k,{h_0}_{|\E^k}))$. We
obtain a period domain $D=G/V$ in the flag manifold
$\check{D}=G_\CC/Q$. Summarizing, considering the classifying map for the
\cvhs\- $E$ on $X$ we obtain~:
\begin{itemize}
\item a real form $\G$ of Hodge type of the Zariski-closure $\GC$ of
  $\rho(\Ga)$ in $\GL(n, \CC)$ such that $\rho(\Gamma) \subset G =
  \G(\RR)$. 
\item a period domain $G/V$ for $G$ and a holomorphic $\rho$-equivariant
  period map $f_\rho: \ti{X} \lo G/V$ making the following diagram
  commutative~:
\begin{equation*}
\xymatrix{
\ti{X} \ar[r]^{f_\rho} \ar[dr]_{\overline{f_{\rho}}} &G/V \ar[d]^\pi \\
&G/K
}
\;\;.
\end{equation*}
\end{itemize}

We can assume without loss of generality that the group $\G$ is $\RR$-simple.
For a simple real Lie group $G$ of Hodge type the continuous cohomology
$H^\bullet_{ct}(G, \RR) $ concentrates in even degree. Cohomology in
degree $2$ is not useful for a general group of Hodge type~: it is
non-zero if and only if $G$ is of Hermitian type (in which case 
the $G$-invariant K\"ahler form $\omega_{G/K}$ generates $H^2_{ct}(G,
\RR)$). On the other hand the cohomology $H^4_{ct}(G, \RR) $ is always
non-trivial for $G$ of Hodge type~: it
contains the first Pontryagin form $p_1(\FF) \in H^4(\lieg, K,
\RR)$ of the real oriented $\mathcal{C}^\infty$ vector bundle $$ \FF := G \times_{K, \Ad} \liek$$ on $G/K$ associated to the
$K$-principal bundle $K \hookrightarrow G \lo G/K$, the adjoint
representation $\Ad : K \lo \Aut(\liek)$. Here the form $p_1(\FF)$ is
computed with respect to the standard $G$-invariant Riemannian metric on $\FF$ 
coming from the Killing form on $\liek$. The main idea for proving theorem~\ref{main} consists in 
studying the pull-back of this class under $\overline{f_{\rho}}$ using
the factorization through the holomorphic horizontal map $f_\rho$. We prove the following~:

\begin{theor} \label{pont}
Let $\rho: \Ga= \pi_1(X) \lo G$ be the monodromy of a complex
variation of Hodge structure $f_\rho: \tilde{X} \lo G/V$ on a compact K\"ahler manifold
$X$. Suppose that $\rho(\Ga)$ is unbounded in $G$. 

Then the $(2,2)$-form
$\overline{f_\rho}^* \, p_1(\FF) \in A^{2,2}(X, \RR)$ is non-positive
on $X$. 

In particular the class $\rho^* p_1(\FF) \in H^4(\Ga, \RR)$ is non-zero
except if $\overline{f_\rho}^*\, p_1(\FF)$ vanishes identically on $X$.
\end{theor}

Some comments on theorem~\ref{pont}~:
\begin{itemize}
\item As far as I know the only non-trivial case of theorem~\ref{pont}
  which was previously deals with $G = Sp(n,1)$. In this case $G/V$ is the twistor space of the quaternionic symmetric
  space $G/K$ and the quaternionic form $p_1(\FF)$ pulls back via $\pi$ to
  $-c_1(L)^2$, where $L$ is a holomorphic $G$-equivariant Hermitian
  line bundle on $G/V$. Moreover $L$ is seen to be positively
  curved in the horizontal directions \cite{sal0}. Thus in this particular
  case the form $\pi^* p_1(\FF)$ is already non-positive on pairs of
  horizontal vectors.
\item for the general case however theorem~\ref{pont} uses the full-strength of the definition of a
  \cvhs~: the form $\pi^* p_1(\FF)$ is shown to be non-positive only on pairs of
  {\em commuting} vectors in the horizontal directions. Thus one
  crucially uses the classical fact that the image of $df_{\rho}(x)$ for any
  point $x \in \tilde{X}$ can be canonically identified with an {\em abelian horizontal} subalgebra of
  $\gc$. 

\end{itemize}

To deduce theorem~\ref{main} from theorem~\ref{pont} it is enough to
understand what happens in the case where the $(2,2)$-form $\overline{f_\rho}^* \,p_1(\FF) \in A^{2,2}(X,
\RR)$ vanishes identically on $X$. This case looks very peculiar. It happens of
course if the map $f_\rho$ factorizes through a curve but it seems
there are other cases which deserve further study. 
However for our purpose it is enough to argue as follows. Given $E=\oplus_{i=0}^k E^i$ our
\cvhs\- on $X$ let $Q_{k} \supset Q$ be the
parabolic subgroup of $G_\CC$ stabilizing $\E_k = \mathbb{F}_k$. Let
$q_k: \check{D} = G_\CC/Q  \lo \check{D}_{k}:= G_\CC /Q_{k}$ be the corresponding
holomorphic projection of flag varieties. By replacing $E$ by $\Lambda^{r_{k}}
E$ (where $r_k = \rk E^k$) we can assume that $r_k=1$. In this case we obtain the following~:

\begin{corol} \label{factor}
 Let $\rho: \Ga= \pi_1(X) \lo G$ be the monodromy of a complex
variation of Hodge structure $E=\oplus_{i=0}^k E^i$ on a compact
K\"ahler manifold, with associated period map $f_\rho: \tilde{X} \lo
G/V$. Suppose that $r_k = \rk E^k =1$. Then~:
\begin{itemize}
\item either the $(2,2)$-form $\overline{f_\rho}^* \, p_1(\FF) \in
  A^{2,2}(X, \RR)$ is non-positive on $X$ and negative at some point
  of $X$. In this case~:
$$\rho^* p_1(\FF) \not = 0\in H^4(\Ga, \RR)\;\;.$$ 
\item or the map $q_k \circ f_\rho: \ti{X} \lo \check{D}_k$ factorizes
  through a complex curve.
\end{itemize}
\end{corol}

In the second case one easily obtains as in \cite{KlKoMau}  that the
pull-back of a well-chosen $G$-invariant real $(1,1)$-form on $\check{D}$ gives a non-trivial
element in $H^2(\Ga, \RR)$, thus finishing the proof of theorem~\ref{main}.

\subsection{Acknowledgements}
It is a pleasure to thank B. Claudon, P.Eyssidieux, V. Koziarz,
J.Maubon, C. Simpson and D.Toledo for their interest
in this work. 

\subsection{Organization of the paper} The paper is organized as
follows. In section~\ref{rappels} we collect some well-known facts
concerning continuous cohomology of real simple Lie groups. In
section~\ref{char} we  justify diagram~(\ref{diag}). In
section~\ref{nonab} we remind the relevant facts from non-Abelian
Hodge theory and in section~\ref{geom} the geometry of period
domains. Section~\ref{first} contains the fundamental computation of
$\pi^* p_1(\FF)$ for pairs of commuting horizontal vectors
(c.f. theorem~\ref{pont1}) and the proof of theorem~\ref{pont}. The
last two sections prove corollary~\ref{factor}, theorem~\ref{main},
theorem~\ref{linear2} and theorem~\ref{lattice}.

\section{Continuous cohomology of real simple Lie groups} \label{rappels}

In this section we collect some classical results for the convenience
of the reader. For simplicity all the cohomologies will be understood
with real coefficients.

\subsection{Cartan's theorem}
Let $G$ be a connected simple real Lie group of non-compact type, $U$
its compact dual and $K$ its maximal compact subgroup.
By Van Est's theorem~:
\begin{equation} \label{e2}
H^\bullet_{ct}(G) \simeq H^\bullet(\lieg, K) \simeq
H^\bullet(U/K) \;\;.
\end{equation}
We are thus reduced to compute the cohomology of $U/K$.
More generally let $V$ be any connected Lie subgroup of $U$. We want
to compute $H^\bullet (U/V)$.

First for $V= \{1\}$. Let $\liet_U$ and $\liet_V$, $\liet_V \subset
\liet_U$, be maximal Abelian
subalgebras of $\lieu$ and $\liev$ respectively. By Hopf theorem
\cite{bor0}~:
$$ H^\bullet(G)= \Lambda^\bullet (z_1, \cdots ,z_n)$$ where the
$z_i$'s are universal transgressive elements. The
Cartan-Chevalley theorem implies~:
$$H^\bullet(BU) = S^\bullet(\liet_U)^{W_{U}}$$ is the ring of
$W_U$-invariants polynomials on $\liet_U$. If $y_1, \cdots , y_n$
denotes the image of $z_1, \cdots, z_n$ by transgression for the
$U$-principal bundle $U \hookrightarrow EU\lo BU$ then the $y_i$'s
generate $H^\bullet(BU)$. 

For any connected $V \subset U$~: let 
$$\rho^*(V, U): H^\bullet(BU)\simeq S^*(\liet_U)^{W_{U}}  \lo
S^*(\liet_V)^{W_{V}} \simeq H^\bullet(BV)$$ be the natural restriction
map. Let $(L,d)$ be the differential algebra defined by~:
$$
L= H^\bullet(BV) \otimes_{\ZZ} H^\bullet(U)$$
and
\begin{equation}
\begin{split}
d(1 \otimes z_i) & = \rho^*(V, U) y_i \otimes 1 \quad (1 \leq i \leq n) \\
d(b \otimes 1) & =0 \quad \text{for} \; b \in H^\bullet(BV) \;\;.
\end{split}
\end{equation}

Let denote by $q: H^\bullet(BV) \lo H^\bullet(L,d)$ be the
homomorphism induced by the inclusion $(H^\bullet(BV),0) \subset
(L,d)$ of $dg$-algebras and by $\sigma^*: H^\bullet(BV) \lo
H^\bullet(U/V)$ the characteristic homomorphism deduced from the
classifying map $U/V \lo BV$.

\begin{theor}[H.Cartan] \cite[p.187]{bor0} \label{cartan}
There exists an isomorphism $\mu : H^\bullet(L, d) \lo H^\bullet(U/V)$
such that the following diagram commutes~:
$$
\xymatrix{
H^\bullet(BV) \ar[r]^q \ar[dr]_{\sigma^*}& H^\bullet(L,d) \ar[d]^\mu \\
& H^\bullet(U/V) \;\;.}
$$
\end{theor}

\subsection{Equal rank case}
\begin{theor} \cite[theor 26.1]{bor0} \label{bor0}
Let $U/V$ be a compact homogeneous space, with $U$  connected and $V$
connected compact. If $\rk U= \rk V$ then~:
\begin{itemize}
\item[(a)] $\rho^*(V, U) : H^\bullet(BU) \lo H^\bullet(BV)$ is injective.
\item[(b)] $H^\bullet(G/U) =
  S^\bullet(\liet_U)^{W_{V}}/<S^\bullet(\liet_U)^{W_{U}}>^+$ where
  $<S^\bullet(\liet_U)^{W_{U}}>^+$ is the ideal of
  $S^\bullet(\liet_U)^{W_{V}}$ generated by elements of
  $S^\bullet(\liet_U)^{W_{U}}$ of positive degrees.
\item[(c)] If $s_i -1$, $1 \leq i \leq n$, resp. $r_j -1$, $1 \leq j \leq
  n$, are the degrees of the $z_i$'s, resp. of universal transgressive
  elements basis (as an exterior algebra) of $H^\bullet(BV)$ then the
  Poincar\'e polynomial of $U/V$ is~:
$$P(U/V, t ) = \frac{(1-t^{s_1}) \cdots (1-t^{s_l})}{(1-t^{r_1})
  \cdots (1-t^{r_l})} \;\; \text{(Hirsch formula)} \;\;.$$
\end{itemize}
\end{theor}

\subsection{Low-degree cohomology}
Consulting tables for the constants $s_i$ and $r_i$ \cite{fomenko} one
easily deduces from theorem~\ref{bor0} the following~: 
\begin{corol}
Let $U/K$ be a symmetric space of compact type with $U$ connected
simple real compact Lie group of dimension larger than $3$. Then~:
\begin{itemize}
\item[(a)] $h^2(U/K) = 0$ except if $U$ is of Hermitian type in which
  case $h^2(U/K)=1$.
\item[(b)] 
$$
h^4(U/K)=
\begin{cases}
0 & \text{if} \; K \; \text{is simple,} \\
2 & \text{if} \; U \; \text{is of Hermitian type,} \\
1 & \text{otherwise}.
\end{cases}
$$ 
\item[(c)] If $\rk K = \rk U$ then $H^\bullet(U/K)$ is concentrated in
  even degrees, $h^2(U/K)=0$ and $h^4(U/K)=1$ except if $U$ is of
  Hermitian type in which case $h^2(U/K)=1$ and $h^4(U/K) =2$.
\end{itemize}
\end{corol}

One can also show by elementary methods~:
\begin{theor}
Let $U/K$ be a symmetric space of compact type with $U$ connected
simple real compact Lie group of dimension larger than $3$ and $\rk U
= \rk K$. Then the first Pontryagin class $[p_1(\FF)] \in H^4(U/K)$ of
the real bundle $\FF= U\times_{K, \Ad} \liek$ does
not vanish.
\end{theor}

\begin{rem}
One also easily shows that if the total Pontryagin class $p(U/K) \in
H^{\bullet}(U/K)$ is non-trivial then necessarily $\rk U = \rk K$ \cite{fomenko}.
\end{rem}

\section{Characteristic classes of representations} \label{char}

Let $\Ga$ and $G$ be topological groups. 
Let $\rho: \Ga \lo G$
be a group morphism. Our goal in this section is to recall some basic
facts concerning the pull-back map it induces in cohomology
$$ \rho^* : H^\bul_c(G, \RR) \lo H^\bul_c(\Ga, \RR)\;\;,$$
in particular the commutative diagram~(\ref{diag}).

Suppose $\Ga$ is a discrete group. Thus $\rho$ factorizes through
$G^\delta$ (the group $G$ with the discrete topology) and the pull-back map in cohomology can be written
$$ H^\bullet_c(G, \RR)  \stackrel{i^*}{\hookrightarrow} H^\bullet(G^\delta, \RR) = H^*(BG^\delta, \RR)
\stackrel{\rho^*}{\lo} H^\bullet(\Ga, \RR) \;\;,$$
where $i: G^\delta \lo G$ is the canonical map and the injectivity of
$i^*$ is a classical result of Borel and Selberg.

We will be interested in the case where $\Ga$ is the fundamental group
of a closed smooth manifold
and $G$ is a real simple Lie group of non-compact type with maximal compact subgroup $K$ and associated
Riemannian symmetric space $G/K$.

\subsection{Continuous cohomology versus discrete cohomology for Lie
  groups} For details on this section we refer to \cite{dupont}.
Let $G$ be as above. The continuous cohomology $H^\bullet_c(G, \RR)$ is
described by Van Est's theorem~:
$$H^\bul_c(G, \RR) \simeq H^\bul(A^*(G/K)^G, d) \simeq H^\bul(\lieg, K,
\RR)$$
where $(A^*(G/K)^G, d)$ denotes the $G$-invariant subcomplex of the
real De Rham complex of $G/K$.

Let $\FF_{G/K} \lo BG^\delta$ denotes the flat bundle with fiber $G/K$
associated to the canonical  flat $G$-bundle $\FF_G$ over $BG^\delta$. As $G/K$ is
contractible the bundle $\FF_{G/K}$ admits a continuous section $s$, unique up to homotopy. Finally one obtains the
following diagram, uniquely defined in the homotopy category~:
\begin{equation} \label{diag1}
\xymatrix{
\FF_{G/K} \ar[d]^{p} \\
BG^\delta \ar@/^1pc/[u]^{s}}
\end{equation}
Equivalently the section $s$ can be interpreted as a
$G^\delta$-equivariant map
$$ EG^\delta \stackrel{s}{\lo} G/K\;\;.$$

Any $G$-invariant closed form $\lambda \in A^i(G/K)$ defines by flat
glueing a closed simplicial differential form $\lambda \in
A^i(\FF_{G/K})$ with closed pull-back $s^*\lambda \in
A^i(BG^\delta)$. Equivalently $s^*\lambda$ is a $G^\delta$-invariant
closed form on $EG^\delta$.

Finally the map
$$H^\bul_c(G, \RR)  \stackrel{i^*}{\lo} H^\bullet(G^\delta, \RR) =
H^*(BG^\delta, \RR) $$
is canonically identified with~:
\begin{equation} \label{e1}
s^* : H^\bul(A^*(G/K)^G, d) \simeq H^\bul(\lieg, K,
\RR) \lo H^\bul(G^\delta, \RR) \;\;.
\end{equation}

\subsection{The characteristic homomorphism}
Let $\Ga$ be a discrete group and $G$ as above. One has the pull-back
diagram, uniquely defined in the homotopy category~:
\begin{equation} \label{diag2}
\xymatrix{
\FF_{G/K, \rho} \ar[d]^{p} \ar[r] & \FF_{G/K} \ar[d]^{p} \\
B\Gamma \ar@/^1pc/[u]^{s_\rho} \ar[r]_{\rho} & BG^\delta \ar@/^1pc/[u]^{s} }
\end{equation}
Equivalently the section $s_{\rho}$ can be interpreted as a
$\rho$-equivariant continuous map
$$ E\Gamma \stackrel{s_{\rho}}{\lo} G/K\;\;.$$

Suppose moreover that $\Ga = \pi_1(M)$ is the fundamental group of a
connected topological space $M$. Gluing the diagram~\ref{diag2} along the
classifying map (canonical in the homotopy category) $c: M \lo B\Ga$ defined by the
universal cover $\ti{M}$ of $M$, one obtains~:
\begin{equation} \label{diag3}
\xymatrix{
c^*(\FF_{G/K, \rho}) \ar[d]^{p} \ar[r] & \FF_{G/K, \rho} \ar[d]^{p} \ar[r] & \FF_{G/K} \ar[d]^{p} \\
M \ar@/^1pc/[u]^{f_\rho}   \ar[r]_c &B\Gamma \ar@/^1pc/[u]^{s_\rho} \ar[r]_{\rho} & BG^\delta \ar@/^1pc/[u]^{s} }
\end{equation}
Equivalently the section $f_{\rho}$ can be interpreted as a
$\rho$-equivariant continuous map
$$ f_{\rho} : \ti{M} \stackrel{c}{\lo}  E\Gamma \stackrel{s_{\rho}}{\lo} G/K\;\;.$$

Finally one has the following commutative diagram, which refines diagram~(\ref{diag})~:
\begin{equation} \label{diag4}
\xymatrix{
H^\bullet_c(G, \RR) \ar@{=}[d] \ar[r]^{i^{*}} & H^\bullet(G^\delta, \RR)
\ar@{=}[d] \ar[r]^{\rho^*} & H^\bullet(\Ga, \RR)  \ar@{=}[d] \ar[r]^{c^{*}} &
H^\bullet(M, \RR) \ar@{=}[d]  \\
H^\bul(A^*(G/K)^G, d) \ar[r]^{s^{*}} \ar@/_1pc/[rr]_{s_\rho^* } \ar@/_3pc/[rrr]_{f_\rho^* } & H^\bullet(A^*(G^\delta, \RR),d) \ar[r]^{\rho^*} & H^\bullet(A^*(B\Ga, \RR),d) \ar[r]^{c^{*}} &
H^\bullet(A^*(M, \RR),d) }
\end{equation}

\section{Non-Abelian Hodge theory} \label{nonab}

In this section $\GC$ is a reductive $\CC$-algebraic group and $G_\CC=
\GC(\CC)$ its Lie group of complex points. We fix once for all a
maximal compact subgroup $U$ of $G_\CC$.
Let $(X, \omega)$ be a compact connected K\"ahler manifold with
fundamental group $\Ga = \pi_1(X)$ (the role of the base point will be 
unimportant for our discussion). 

\subsection{Simpson's correspondence} 
Non-Abelian Hodge theory establishes a correspondence between (reductive)
$G_\CC$-principal local systems on $X$ (or equivalently $G_\CC$-conjugacy
classes of representations $\Ga \lo G_\CC$ and (semi-harmonic) $\GC$-principal Higgs bundles
on $X$~:

\begin{defi}
A $\GC$-principal Higgs bundle on $X$ is a pair
$(\holP, \theta)$, where 
\begin{itemize}
\item $\holP$ is a principal holomorphic $G_\CC$-bundle on $X$.
\item $\theta \in \Ad \holP \otimes \Omega^1_X$ satisfies
  $[\theta, \theta] =0$ (where $\Ad \holP := \holP
  \times_{G_{\CC}} \lieg_\CC$).
\end{itemize}
A $\GC$-principal Higgs bundle $(\holP, \theta)$ on $X$ is said
to be of semi-harmonic type if its Chern classes vanish and for some
irreducible $\GC$-module $V$ (and then for any) the Higgs vector-bundle
$\holP \otimes_{G_{\CC}} V$ is Higgs semi-stable.
\end{defi}

The correspondence is obtained as follows.
Let $\rho: \Ga \lo G_\CC$ be a reductive representation and $\PP=(P, D)$
the associated flat complex $G_\CC$-bundle on $X$. Here $P$ denotes the
principal $G_\CC$-bundle $\tilde{X} \times_{\rho}G_\CC$ with flat 
connection $D \in A^1(\Ad P)$. As $X$ is compact K\"ahler and $\rho$ is reductive there exists an
essentially unique $\rho$-equivariant harmonic map $$f: \tilde{X} \lo 
G_\CC/U$$ (where $G_\CC/U$ denotes the symmetric space of $G_\CC$),
defining an {\em harmonic $U$-reduction} $P_U$ of $P$. Decompose the
flat connexion $D$ as
$$D = \nabla + \alpha\;\;,$$ 
where $\nabla$ is the canonical connexion on the $U$-principal bundle $P_U$ and $\alpha \in
A^1(X, \Ad P)$. Decompose furthermore using types~:
\begin{equation*}
\begin{split}
\nabla &= \partial_K + \delbar \;\;, \\
\alpha &= \theta + \theta^*\;\;,
\end{split}
\end{equation*}
where $\partial_K$ is of type
$(1,0)$, $\overline{\partial}$ is of type $(0,1)$, $\theta \in A^{1,0}(\Ad P)$
and $\theta^* = \tau(\theta) \in A^{0,1}(\Ad P)$ is the conjugate of
$\theta$ with respect to the $U$-reduction.
Define $D' = \partial_K + \theta^*$, $D'' = \delbar + \theta$, thus $D
= D' + D''$. As $D$ is flat and the $U$-reduction $P_U$ is harmonic,
$(D'')^2=0$, that is~: 
$$\delbar^2 =\delbar(\theta) = [\theta, \theta]=0 \;\;.$$
Finally $(\holP = (P, \delbar) , \theta)$ is a
$\GC$-Higgs bundle (of semi-harmonic type). Notice that knowing
$(\holP, \theta)$ is equivalent to knowing $(P, D'')$.

\sspace
Let $\GC-dR$ be the differential graded-category of
flat $G_{\CC}$-bundles on $X$~: an object is a flat bundle
$\PP=(P, D)$ on $X$ and $$\HHom_{\GC-dR}(\PP, \PP') =
(A^\bullet(\Hom(\Ad P, \Ad P'), D_{\Hom(P,P')})\;\;.$$ 
Let $\GC-Dol$ be the differential graded-category of
semi-harmonic $\GC$-Higgs bundles on $X$~: an object is a semi-harmonic
$\GC$-Higgs bundle $(\holP, \theta_{\holP})$ on $X$ and
$$\HHom_{\GC-Dol}(\holP, \theta_{\textnormal{P}}), (\holP', \theta_{\holP'})) =
(A^\bullet(\Hom(\Ad P, \Ad \mathcal{P'}), D''_{\Hom(P,P')})\;\;.$$
\begin{theor}[Simpson] \label{equi}
The functor $F: \GC-dR \lo \GC-Dol$ associating to the flat bundle $\PP$
the Higgs bundle $(\holP, \theta)$ is a quasi-equivalence of
differential graded categories.
\end{theor}

It implies the geometric (weaker) version \cite[theor. 9.11 and lemma 9.14]{si2}~:
\begin{theor}[Simpson] \label{sim}
The functor $F$ induces a real-analytic diffeomorphism
$$ \phi_\G:\M(\Ga, \GC)(\CC)  \lo \MD(X, \GC)(\CC) $$
between the Betti moduli space $\M(\Ga, \G)(\CC) = (\HHom(\Ga,
\GC)//\GC^\ad)(\CC)$ of representations of $\Ga = \pi_1(X)$ in $G_\CC$
and the Dolbeault moduli space of $\GC$-Higgs bundles of semi-harmonic type.
\end{theor}

\subsection{$\G$-variations of Hodge structures}  \label{Gcvhs}
The moduli space $\MD(X, \G)$ of $\GC$-Higgs bundles of semi-harmonic
type carries a natural $\CC^*$-action~: an
element $t \in \CC^*$ maps $[(\holP, \theta)]$ to
$[(\holP, t \cdot \theta)]$ \cite[p.62]{si2}.
The fixed points of this action are of particular
importance~: they are systems of $\G$-Hodge bundles \cite[p.44]{si1}
and correspond by Simpson's correspondence to (isomorphism classes of)
$\G$-complex variations of Hodge structure ($\G$-\cvhs).

\subsubsection{Hodge datum} \label{Hodge_datum}
\begin{defi} 
A (pointed) Hodge datum is a pair $(\G, u)$, where $\G$
is a real reductive algebraic group and $u: \U(1) \lo \Aut(\G)^0
\subset \G^\ad$ is a
morphism of real algebraic groups such that $C= u(-1)$ is a Cartan
involution of $\G$ (that is~: $C^2 =1$ and $\tau := C \sigma = \sigma
C$ is the conjugation of $\G_\CC$ with respect to a compact real form $\U$,
where $\sigma$ denotes the conjugation of $\G_\CC$ with respect to $\G$).
\end{defi}

In particular the Cartan involution $C$ is inner. One
easily shows (c.f. \cite[section 4.4]{si1}) that an algebraically connected real reductive group
$\G$ admits a Hodge datum $(\G, u)$ (one says that $\G$ is of {\em
  Hodge type}) if and only if $\G$ contains an anisotropic maximal
torus $\TT$. In other words the reductive real Lie group $\G(\RR)$
has the same real rank than any of its maximal compact
subgroups.

\subsubsection{Period domains} \label{period}

\begin{defi}
Let $(\G, u )$ be a Hodge datum. We denote by~:
\begin{itemize}

\item $V= Z_{G}(u)$ the centralizer of $u$ in $G= \G(\RR)$. As $V$ is invariant
by $C$, $V$ is contained in $U = \U(\RR)$, in particular $V$ is
compact.
\item $K$ the centralizer of $C$ in $L$. Notice that $K$ co\"incide with the
  intersection $U \cap L$. Thus the group $K$ is a maximal compact
  subgroup of $L$. 
\end{itemize}
\end{defi}

Finally the Hodge datum $(\G, u)$ defines canonically the chain of
inclusions of compact groups $V \subset K \subset U$.

\sspace
Let $(\G, u)$ be a Hodge datum. Let $\lambda : \G \lo \GL(E)$ be a
real (resp. complex) representation $\lambda : \G \lo \GL(E)$ of 
$\G$. If $\G= \G^\ad$ (or more generally if $\lambda$ factorizes
through $\G^\ad$) the composite $\lambda \circ u : \U(1)\lo \G^\ad
\lo \GL(E)$ defines a weight $0$ real (resp. complex) Hodge structure on
$E$ polarized by $\lambda \circ u (-1)$.
In particular the adjoint representation of $\G$ defines on the Lie
algebra $\lieg$ a weight $0$ polarized real Hodge structure~:
$$\lieg_\CC = \bigoplus_{i \in \ZZ} \lieg_\CC^{i, -i} \;\;,$$
where $u(z)$ acts on $\lieg_\CC^{i,-i}$ via multiplication by $z^{-i}$.
We will denote by $F^\bullet \lieg_\CC$ the corresponding decreasing
Hodge filtration. The polarization is given 
by the Killing form $\beta_\G$.

\begin{defi}
One denotes by $\qc \subset \lieg_\CC$ the Lie sub-algebra $F^0 \lieg_\CC$
and $Q \subset G_\CC$ the corresponding subgroup.
\end{defi}
One easily check that $\qc$ is a parabolic sub-algebra of $\lieg_\CC$,
with Levi sub-algebra $\vc$ the complexified Lie algebra of $\liev = \text{Lie}(V)$.

\begin{defi}
Let $(\G,u)$ be a Hodge datum. The period domain $D$ associated to $(\G,
u)$ is the $G$-conjugacy class of $u$.
\end{defi}

Thus $D$ naturally identifies with $G/V$. Let $\check{D} = G_\CC/Q$ be
the flag manifold of $G_\CC= \G(\CC)$ defined by $Q$, the natural morphism
$D = G/V \hookrightarrow \check{D} = G_\CC/Q$ is an open embedding and thus defines
a natural $G$-invariant complex structure on $D$.

\subsubsection{Horizontality}
The holomorphic tangent bundle $TD$ naturally identifies with the
$G$-equivariant bundle $(G_\CC \times_Q \lieg_\CC/\qc)_{|D}$.
\begin{defi}
The superhorizontal tangent bundle $T_{sh}D$ is the holomorphic sub-bundle
$(G_\CC \times_Q F^{-1}\lieg_\CC/\qc)_{|D}$ of $TD$.
\end{defi}

\begin{defi} \label{hor}
Let $\mathcal{C}$ be the category whose objects are pairs $(Y, R_Y)$, where $Y$ is a
complex smooth analytic space, $R_Y \subset TY$ a holomorphic distribution,
and a morphism $f: (Y, R_Y) \lo (X, R_X)$ in $\mathcal{C}$ is a holomorphic {\em
  horizontal} map $f: X \lo Y$~: one requires that $df(R_Y) \subset
R_X$. We will look at the category of smooth analytic spaces 
as a subcategory of $\mathcal{C}$, the distribution being the full
tangent space.
\end{defi}

\subsubsection{$\GC$-\cvhs}
With all these definitions we can define the main actors in Simpson's
theory~:
\begin{defi}
Let $X$ be a complex analytic manifold with fundamental group
$\Ga$ and universal cover $\tilde{X}$. Let $\G_\CC$ be a complex
reductive algebraic group.
A $\G_\CC$-complex variation of Hodge structure ($\G_\CC$-\cvhs) is a 
Hodge datum $(\G,u)$ with period domain $D$ for a real form $\G$ of
$\GC$, a representation $\rho: \Ga \lo
G=\G(\RR) \subset G$ (called the monodromy of
the variation) and a holomorphic {\em
  horizontal} $\rho$-equivariant map $f: 
\tilde{X} \lo D$ (called {\em period map}).
\end{defi}

\begin{defi} \label{bundles}
Let $(\G, u, \rho, f: \tilde{X} \lo D)$ be a
$\G_\CC$-\cvhs\- on $X$. Let $\tau: Q \lo V_\CC
\stackrel{i}{\hookrightarrow} G_\CC$ be the reduction of $i$ to the Levi
$V_\CC$ of $Q$. One associates to it the following principal
$G_\CC$-bundles on $X$~:

\begin{itemize}
\item the flat $G_\CC$-bundle $\PP = \tilde{X} \times_{\Ga,
    \rho} G_\CC$, which is naturally a holomorphic bundle. Notice that
  this holomorphic structure is compatible with the identification $\PP:= f^*((L_\CC
  \times_{Q,i} G_\CC)_{|D})$ (descent to $X$ of the) pull-back via $f$ of the holomorphic 
  $G_\CC$-bundle $G_\CC \times_{Q,i} G_\CC$ on $\check{D}$.

\item the holomorphic $G_\CC$-bundle $\holP= f^*((L_\CC \times_{Q, \tau}
  G_\CC)_{|D}$ (descent to $X$ of the) pull-back via $f$ of the holomorphic 
  $G_\CC$-bundle $L_\CC \times_{Q,\tau} G$ on $\check{D}$.
\end{itemize}
\end{defi}

\subsection{$\G$-\cvhs\ and Simpson correspondence}

Let $(\G, u: \U(1) \lo \G^\ad, \rho:
\Ga \lo G, f: \tilde{X} \lo D)$ be a
$\G_\CC$-\cvhs\- on $X$.

On the Betti side it canonically defines the isomorphism class of the flat $G$-bundle
$\PP$.

On the Dolbeault side~: the adjoint bundle
$\Ad \holP = \holP \times_{G_\CC, \Ad} \lieg_\CC= f^*((G_\CC \times_{Q, \Ad
   \circ \tau}\lieg_{\CC})_{|D})$ identifies with the graded bundle
$\Gr_F \,\Ad \PP$ of the weight $0$ complex variation of Hodge structure
 $\Ad \PP = \mathcal{E}_{\lieg}$ associated to the 
 representation $\G_\CC \lo \GL(\lieg_\CC)$~:
$$\Ad \holP = \bigoplus_{p \in \ZZ} (\Ad \holP)^{p, -p}\;\;.$$

\begin{defi}
Define $\theta_f \in (\Ad \holP)^{-1, 1} \otimes \Omega^1_X$ as the
differential of $f$.
\end{defi}

Thus $(\holP , \theta)$ is a semi-stable $\GC$-Higgs-bundle and the
$\GC$-\cvhs\- $(i: \G \hookrightarrow \Res \G, \rho, f: 
\tilde{X} \lo D)$ canonically defines the point $[(\holP, \theta)]$
(called a system of $\G$-Hodge 
bundles by Simpson) in $\MD(X,\G)$.

\begin{prop}\cite[cor 4.2]{si1} \label{simpson} \label{prop1}
Let $[\rho] \in \M(\Ga, \GC)(\CC)$ with $\rho$ reductive. Then $\phi_{\GC}([\rho]) \in \MD(X, \GC)(\CC)$ is
$\CC^*$-fixed if and only if $\rho$ is the monodromy of a
$\GC$-complex variation of Hodge structure $f: \tilde{X} \lo D$.
Moreover $\phi_{\GC}([\rho]) = [(\holP, \theta)]$.
\end{prop}

\section{Geometry of period domains} \label{geom}

\subsection{The fibration $D=G/V\lo G/K$}
The action of the compact form $U \subset G_\CC$ on $\check{D}$ is
transitive with stabilizer $V$ at the point $u \in
\check{D}$. The inclusion $V \subset K \subset G$ (resp. $V
\subset K \subset U$) defines a $G$-equivariant
(resp. $U$-equivariant) $\mathcal{C}^{\infty}$
fibration $\pi: D=G/V \lo 
G/K$ (resp. $\check{\pi}: \check{D}= U/V \lo 
U/K$) over the symmetric space $G/K$ (resp. its dual
$U/K$)~:

$$
\xymatrix{
D = G/V \ar@{^{(}->}[r]^i \ar[d]_{\pi} & \check{D} = G_\CC/Q & \check{D} =
U/V \ar@{=}[l] \ar[d]^{\check{\pi}} \\
X= G/K  & & U/K }\;\;.$$

The symmetric space $G/K$ (resp. $U/K)$ does not admit any
$G$-invariant (resp. $U$-invariant)
complex structure except in the case where $G$ is of Hermitian type
(i.e. when the Hodge structure on $\lieg$ is of type $(-1,1)$, $(0,0)$
and $(1,-1)$). Even when $G$ is of Hermitian type the projections $\pi: D \lo G/K$ and $\check{\pi}: \check{D}
  \lo U/K$ are not holomorphic in general. However their fibers (in
  particular the fibre $K/V$ over the base point $u \in D$) are always
  compact holomorphic subvarieties of $D$ and $\check{D}$
  \cite[p.261]{gsc}. In particular $D$ is a Stein manifold only in the
  case where $D=X$ is a Hermitian symmetric space. Let us also
  emphasize that the projections $\pi: D \lo G/K$ and
  $\check{\pi}_{|D}: D \lo U/K$ do not co\"incide (even if they have
  the same fibre $K/V$ at the base point $u \in D$).

\subsection{Root data}
Let $(\G,u)$ be a Hodge datum with Cartan involution $C$. Let $V
\subset K \subset U$ be the compact subgroups of $G_\CC$ it defines
(c.f. section~\ref{Hodge_datum}). We will also denote by $\sigma$ and
$\tau$ the complex conjugations of $\gc$ with respect to $\lieg$ and
$\lieu$ respectively. Let $T\subset V$ be a maximal compact torus of
$G$ containing $u(U(1))$. Let $\liet$ be its Lie algebra and $\hlc =
\liet \otimes \CC$ be the corresponding Cartan subalgebra of
$\gc$. The Cartan subalgebra $\hlc$ is
$\tau$-stable and $\sigma$-stable. 

\lspace
The adjoint representation of $\hlc$ on $\gc$ defines a root space
decomposition 
$$ \gc = \hlc \oplus \bigoplus_{\alpha \in \Delta} \lieg_{\alpha}\;\;,$$
where $\Delta= \Delta(\gc, \hlc) \subset \hlc^*$ denotes the root
system of $\hlc$ in $\ggc$ and each root space  
$$\lieg_{\alpha} = \{x \in \gc \;/\; [h, x] = \alpha(h) x \;
\textnormal{pour tout } h\in \hlc\}$$ is of complex dimension $1$. Let
$\Delta(\vc, \hlc) \subset \Delta(\kc, \hlc) \subset
\Delta$ be the root subsystems associated to $(\vc, \hlc)$ and $(\kc,
\hlc)$ respectively. We denote by $<\cdot, \cdot>$ the natural scalar
product on $\hlc^*$.

\lspace
Let $\lieg = \kg \oplus \liep$ be the Cartan decomposition of $\lieg$
defined by the Cartan involution $C$. The involution $C$ commutes with
the adjoint action of $\hlc$ thus any $\lieg_{\alpha}$, $\alpha \in
\Delta$ is either contained in $\kc$ (the root 
$\alpha$ is said to be {\em compact}) or in $\pc$ (the root $\alpha$
is said {\em non-compact}). We denote by $\Delta_{\kg}\subset \Delta$
(resp. $\Delta_{\liep} \subset \Delta$),
the set of compact roots (resp. non-compact roots).

\lspace
Let $\Delta_+ \subset \Delta$ be the unique set of positive roots such
that
$\lieg_{-\alpha} \subset \q$ for every $\alpha \in \Delta_{+}$. Let
$\Psi = \{\alpha_1 , \cdots ,\alpha_r\}$ be the corresponding set of
simple roots. For every subset $\theta \subset \Psi$ we denote by~:
\begin{itemize}
\item $|\theta|$ the cardinal of $\theta$.
\item $\theta^c:= \Psi \setminus \theta$ and 
  $<\theta>$ the set of roots in $\Delta_+$ linear combination of
  elements of $\theta$.
\item $\n_{\theta, -}$ the nilpotent Lie algebra $\bigoplus_{\alpha \in
    \Delta_{+} \setminus <\theta>} \lieg_{-\alpha}$.
\item $\liel_{\theta}$ the reductive subalgebra centralizer of
  $\cap_{\alpha \in \theta} \ker \alpha$, thus
$$ \liel_{\theta} = \hlc \oplus \bigoplus_{\alpha \in <\theta>}
(\lieg_{-\alpha} \oplus \lieg_{\alpha})\;\;$$
\item $\q_{\theta}$ the parabolic subalgebra $\liel_{\theta}
  \oplus \n_{\theta, -}$ and $\Qtheta
  \subset \G_{\CC}$ the corresponding parabolic subgroup.
\end{itemize}

\begin{rem} We followed the convention that 
$\Qtheta$ is of parabolic rank $|\theta^c|$. Thus $Q_{\Psi}= \ggc$ and
$Q_{\emptyset}=B$ is the standard Borel subgroup of $\ggc$
associated to $\Delta_{+}$.
\end{rem}

For $\theta \subset \Psi$ we denote by $\check{D}_{\theta}$ the flag
variety $\check{D}_{\theta} = \ggc/ Q_{\theta}$. For 
$\theta \subset \theta' \subset \Psi$ let 
$\pi_{\theta \theta'} : \dtheta \lo \check{D}_{\theta'}$ be the
natural holomorphic projection. 

\lspace
As the subalgebra $\vc$ is reductive and $\hlc \subset \vc \subset
\kc$ there exists a unique subset $\Phi \subset \Delta_k \cap \Psi$
such that
$$ \vc = \liel_{\Phi} \qquad \textnormal{and} \qquad \q =
\q_{\Phi}\;\;.$$ 

\lspace
The group $U$ is a compact form of $G_\CC$ thus the roots $\alpha\in
\Delta$ are real on $\lieh_{\RR}:= \sqrt{-1}\liet \subset
\hlc$. This identifies $\Delta$ with a subset of the dual space
$\lieh_{\RR}^*$. The root system $\Delta$ cuts 
$\lieh_{\RR}$ (and also using the Killing form the dual space
$\lieh_{\RR}^*$) in a disjoint union of facets. The open facets are
the Weyl chambers. Let $\lieh_{\RR, +}$ be the Weyl chamber of
$\lieh_\RR$ defined by $\Psi$. Let $E \in \liet$ be the generator of
$u: U(1) \lo G$. The subset $\Phi \in \Psi$ naturally identifies with
the set of roots $\alpha$ satisfying $\alpha(\sqrt{-1} E)
\geq 0$. 

\lspace
Let $W$ be the Weyl group of $\Delta$, $s_\alpha \in W$ be the
reflection associated to $\alpha$ and $w_0$ be the element of $W$ of
maximal length.

\subsection{Automorphic bundles}

\begin{defi}
Let $\pi: Q_{\Phi} \lo GL(E)$ be a finite dimensional complex
representation. We denote by $\check{\mathcal{F}}_{\pi} := \ggc \times_{Q_{\theta}, \pi}
E$ the holomorphic $\ggc$-equivariant bundle on $\check{D}$ with fibre $E$
associated to $\pi$ and by  $\FF_\pi$ its restriction to $D$.
\end{defi}

Let $\nmoins= \n_{\Phi, -}$ be the nilpotent radical of the
  parabolic subalgebra $\q= \q_{\Phi}$. Then $$\nplus := \sigma(\nmoins) =
  \tau(\nmoins) = F^{-\infty} \lieg_\CC / F^0 \lieg_\CC = \oplus_{\alpha
    \in \Delta_+ \setminus <\Phi>} \lieg_{\alpha}\;\; .$$ 
The complex vector space $\nplus$ is a natural $\q$-module and
$T\check{D}= \check{\FF}_{\nplus}$, $TD= \FF_{\nplus}$. 

\begin{defi}
We will denote by $\Delta^{-1,1} \subset \Delta_+$ the unique set of
positive roots such that 
$$\gc^{-1,1} = \bigoplus_{\alpha \in \Delta^{-1,1}} \lieg_{\alpha}
\;\;.$$
For convenience we will put a total order $<_t$ on $\Delta^{-1,1}$ (we
will distinguish it from the partial order $\leq$ on
$\Delta$ defined by $\Psi$).
\end{defi}

\begin{rem}
The set $\Delta^{-1,1} $ admits a natural decomposition corresponding
to a decomposition of the sub-$\q$-module $\gc^{-1,1}=F^{-1}\gc /F^{0}
\gc $ into a sum of blocs~: 
$$\gc^{-1,1} = \bigoplus_{\beta \in \Phi^c} \lieg^{-1,1}_{\CC,
  \beta}\;\;,$$
with
$$\lieg^{-1,1}_{\CC, \beta} :=\bigoplus_{\gamma \in( \beta +
  <\theta>) \cap \Delta_+} \lieg_ \gamma\;\;.$$ 
Correspondingly~:
$$T_{sh}D = \FF_{\gc^{-1,1}}= \bigoplus_{\beta \in \Phi^c}
\FF_{\lieg^{-1,1}_{\CC, \beta}}\;\;,$$ which corresponds to the usual Hodge
  bloc decomposition $\bigoplus \HHom(H^{p,q}, H^{p-1, q+1})$.
\end{rem}

\subsection{Curvature of automorphic bundles}  \-\- \label{courbure}
As a complex bundle the holomorphic fiber bundle $\FF_\pi =
\check{\FF_{\pi}}$ can be written $G\times_{V, \pi_{|V}} E$. As the
group $V$ is compact on can put on $E$ a $V$-invariant Hermitian
metric which induces a $G$-invariant Hermitian metric on
$\FF_{\pi}$. The Chern connexion on the Hermitian holomorphic bundle
$\FF_{\pi}$ is $G$-invariant and is deduced from the connexion form of
the $V$-principal bundle $V \hookrightarrow G \lo G/V=D$. We will
denote by $\Theta_D(\pi) \in A^{1,1}(D) \otimes \End(\FF_{\pi})$ its
$G$-invariant curvature form. These curvature forms are computed in 
\cite{gsc}.

\begin{lemma} \cite[p.265]{gsc} \label{baseweyl}
There exists a basis $(e_{\alpha} \in \lieg_{\alpha}, h_{\alpha} \in
\lieh_{\RR})_{\alpha \in \Delta}$ of $\gc$ (called a Weyl basis)
satisfying the following properties~:
\begin{itemize}
\item[(a)] $(e_\alpha, e_\beta) = \delta_{\alpha, -\beta}$,
  $[e_{\alpha}, e_{-\alpha}] = h_\alpha$.
\item[(b)] $(h_{\alpha}, x) = \alpha(x)$ for $x \in \hlc$.
\item[(c)] $[e_\alpha, e_\beta]=0 $ if $\alpha \not = - \beta$ et
  $\alpha + \beta \not \in \Delta$.
\item[(d)] $[e_\alpha, e_\beta]=  N_{\alpha, \beta} \, e_{\alpha +
    \beta}$ si $\alpha, \beta, \alpha + \beta \in \Delta$, where
  $N_{\alpha, \beta} \in \RR$, $N_{-\alpha, - \beta}= -N_{\alpha,
    \beta}$.
\item[(e)] $\tau(e_{\alpha})= - e_{\alpha}$.
\item[(f)] $\sigma(\alpha) = \varepsilon_\alpha e_{-\alpha}$ where
  $\varepsilon_\alpha = -1$ if $\alpha$ is compact and 
  $\varepsilon_\alpha = 1$ if $\alpha$ is non-compact.
\item[(g)] $\varepsilon_{\alpha + \beta}= - \varepsilon_\alpha
  \varepsilon_\beta$ if $\alpha, \beta, \alpha+ \beta \in \Delta$.
\end{itemize}
\end{lemma}

\begin{defi}
Let $(e_{\alpha} \in \lieg_{\alpha}, h_{\alpha} \in
\lieh_{\RR})_{\alpha \in \Delta}$ be a Weyl basis.
We denote by $w_\alpha \in \lieg_{\alpha}^*$, $\alpha \in \Delta$,
the linear form dual to $e_{\alpha}$. 
\end{defi}

\begin{theor} \cite[p.269]{gsc} \label{griffiths}
\begin{itemize}
\item[(1)]
The curvature form $\Theta_D$ for the natural connection on
the $V$-principal bundle $V \hookrightarrow G \lo G/V=D$ is
$$\Theta_D = - \sum_{\alpha, \beta \in \Delta^+ \setminus
  \Delta_{\kc}} [e_\alpha, e_{-\beta}]_{\vc} \otimes w_{\alpha}
\wedge \overline{w_{\beta}} + \sum_{\alpha, \beta \in (\Delta(\kc) \cap \Delta_+) \setminus
    <\Phi>}  [e_\alpha, e_{-\beta}]_{\vc} \otimes w_{\alpha}
\wedge \overline{w_{\beta}} \;\;.$$
\item[(2)] Let $\pi: Q_{\Phi} \lo GL(E)$ be a finite dimensional complex
representation. The curvature form $\Theta_D(\pi) \in A^{1,1}(D)
\otimes \End \FF_{\pi}$ for the natural Hermitian
connection on the automorphic vector bundle $\FF_\pi$ is 
$$\Theta_D(\pi) = - \sum_{\alpha, \beta \in \Delta^+ \setminus
  \Delta_{\kc}} \pi([e_\alpha, e_{-\beta}]_{\vc} )\otimes w_{\alpha}
\wedge \overline{w_{\beta}} + \sum_{\alpha, \beta \in (\Delta(\kc) \cap \Delta_+) \setminus
    <\Phi>}  \pi([e_\alpha, e_{-\beta}]_{\vc}) \otimes w_{\alpha}
\wedge \overline{w_{\beta}} \;\;.$$
\end{itemize}
\end{theor}

\begin{corol} \label{courb}
Let $\pi: Q_{\Phi} \lo GL(E)$ be a finite dimensional complex
representation. Let $\xi \in T_{sh}D \simeq \lieg_\CC^{-1,1}$ be a
superhorizontal holomorphic vector
field. Then 
$$\Theta_D(\pi)(\xi \wedge \overline{\xi}) = - [\xi, \sigma(\xi)] \in
\vc= \lieg_\CC^{0,0} \;\;.$$
\end{corol}

\section{First Pontryagin class~: proof of theorem~\ref{pont}} \label{first}

\subsection{Computation of $\pi^*p_1(\FF)$}
Let $\pi: D= G/V \lo G/K$ be the canonical projection. Let
$$ \FF := G \times_{K, \Ad} \liek$$ be the real oriented
$\mathcal{C}^\infty$ vector bundle on $G/K$ associated to the
$K$-principal bundle $K \hookrightarrow G \lo G/K$ and the adjoint
representation $\Ad : K \lo \Aut(\liek)$. We denote by $p_1(\FF) \in
H^4(\lieg, K, \RR)$ the first Pontryagin form of the bundle $\FF$ with
its canonical metric associated to the $K$-invariant scalar product
$-B_\liek(X,Y)$. By definition $p_1(\FF) = -c_2(\FF \otimes_{\RR}
  \CC)$ the second Chern form of the complexified bundle $\FF \otimes_{\RR}
  \CC$. Notice that the pull-back complex bundle $\pi^*\FF \otimes_{\RR}
  \CC$ with its pulled-back Hermitian metric is in fact an automorphic bundle on $G/V$~:
$$\pi^* \FF \otimes_{\RR} \CC  = \FF_{\liek_{\CC}} \;\;.$$
Thus $$\pi^* p_1(\FF) = p_1(\FF_{\liek}) = -c_2(\FF_{\liek_{\CC}} ) \in H^4(\lieg, V, \RR)$$ lies in $H^{2,2}(\lieg,
V, \RR)$.

\begin{theor} \label{pont1}
Let $\xi, \eta \in \gc^{-1,1}$ satisfying $[\xi, \eta]=0$. Then
\begin{equation*}
- 8 \pi^2 \cdot \pi^* p_1(\FF) (\xi \wedge \eta \wedge \ov{\xi} \wedge
\ov{\eta}) = \sum_{ \substack{\alpha, \beta \in \Delta^{-1, 1} \\ \alpha <_t \beta}} (N_{\alpha, - \beta}^2 + <\alpha, \beta> +
N_{\alpha, \beta}^2)  \cdot |\xi^\alpha \eta^\beta - \xi^\beta
\wedge \eta^\alpha|^2 \qquad.
\end{equation*}

\end{theor}

\begin{proof}
As usual~:
$$\pi^*p_1(\FF) = -c_2(\FF_{\kc}) = - \frac{1}{8 \pi^2} \Tr
\Theta_D(\kc)^2 \;\;.$$
For simplicity we will denote by $\Xi$ the $(2,2)$-form $-8 \pi^2 \cdot \pi^*p_1(\FF)$. 
As recalled in corollary~\ref{courb} for any $\xi \in \gc^{-1,1}$ one
has 
$$\Theta_D(\kc)(\xi \wedge \overline{\xi}) = - \ad_{\kc} [\xi,
\sigma(\xi)] \in \End \FF_{\kc}\;\;.$$
Notice that $\sigma = -\tau$ on $\gc^{-1,1}$. Thus for any $\xi, \eta
\in \gc^{-1,1}$ we get~:
\begin{equation}
\Xi (\xi \wedge \eta \wedge \overline{\xi} \wedge \overline{\eta})
= \Tr (\ad_{\kc} [\xi, \tau(\xi)] \ad_{\kc}
[\eta, \tau(\eta)] | \kc) - \Tr(\ad_{\kc} [\xi, \tau(\eta)] \ad_{\kc}
[\eta, \tau(\xi)] | \kc)) \;\;,
\end{equation}
which can be rewritten as~:
\begin{equation}  \label{equa_p_1}
\Xi(\xi \wedge \eta \wedge \overline{\xi} \wedge \overline{\eta})
=  B_{\kc}([\xi, \tau(\xi)], 
[\eta, \tau(\eta)]) - B_{\kc}( [\xi, \tau(\eta)],
[\eta, \tau(\xi)] ) \;\;.
\end{equation}

Let $\alpha, \beta, \gamma, \varepsilon \in \Delta^{-1,1}$ and $\ea,
\eb, \eg, \ee \in \gc^{-1,1}$ the corresponding vectors of our Weyl basis. 
Thus  
\begin{equation*}
\begin{split}
\Xi_{\alpha \beta \ov{\gamma} \ov{\varepsilon}} &=B_{\kc}([\ea,
e_{-\gamma}], [\eb, e_{\varepsilon}]) - B_{\kc}([\ea, e_{-\varepsilon}],
[e_{\beta}, e_{- \gamma}]) \\
& = N_{\alpha, - \gamma} N_{\beta, -\delta} B_{\kc}(e_{\alpha -\gamma}, e_{\beta- \varepsilon}) - N_{\alpha, - \varepsilon}
N_{\beta, - \gamma} B_{\kc}(e_{\alpha -\varepsilon}, e_{\beta -\gamma})
\;\;,
\end{split}
\end{equation*}
with the convention $e_{\alpha - \alpha} = h_{\alpha}$. 

Thus~:
\begin{equation} \label{coeff}
\Xi_{\alpha\beta \ov{\gamma} \ov{\varepsilon}} = 
\begin{cases}
 N_{\alpha, - \gamma} N_{\beta, -\delta} - N_{\alpha, - \varepsilon}
N_{\beta, - \gamma} & \text{if} \; \alpha- \gamma = \beta - \varepsilon
\; \text{and} \; \alpha \not = \gamma,\\
N_{\alpha, - \beta}^2 + <\alpha, \beta> & \text{if} \; \alpha = \gamma
\; \text{and} \; \beta = \varepsilon, \\
0 & \text{otherwise}.
\end{cases}
\end{equation}

This expression can be simplified using the following~:
\begin{lemma} \cite[p.147]{helg}
Let $\alpha$, $\beta$, $\gamma$ and $\varepsilon$ be $4$ roots in
$\Delta$, no two of which have sum $0$. Then
$$ N_{\alpha, \beta} N_{\gamma, \varepsilon} + N_{\beta, \gamma}
N_{\alpha, \varepsilon} + N_{\gamma, \alpha} N_{\beta, \varepsilon} =0
\;\;.$$
\end{lemma}

Let us apply this lemma to the first case of equation~(\ref{coeff})
and the sum of roots $\alpha+ \beta- \gamma - \delta =0$. Using $N_{- \gamma, \alpha}=   - N_{\alpha, - \gamma}$
and $N_{- \gamma, - \varepsilon} = - N_{\gamma, \varepsilon}$ one
obtains~:

\begin{lemma}
\begin{equation} \label{coeff1}
\Xi_{\alpha \beta \ov{\gamma} \ov{\varepsilon}} = 
\begin{cases}
- N_{\alpha, \beta} N_{\gamma, \varepsilon} & \text{if} \; \alpha- \gamma = \beta - \varepsilon
\; \text{and} \; \alpha \not = \gamma,\\
N_{\alpha, - \beta}^2 + <\alpha, \beta> & \text{if} \; \alpha = \gamma
\; \text{and} \; \beta = \varepsilon, \\
0 & \text{otherwise}.
\end{cases}
\end{equation}
\end{lemma}

Let now $\xi, \eta \in \gc^{-1,1}$. We write $\xi = \sum_{\alpha} \xi^\alpha
e_{\alpha}$ and $\eta  = \sum_{\beta} \eta^\beta e_{\beta}$, where for
simplicity of notations all the indices $\alpha, \beta, \gamma,
\varepsilon$ from now on vary in $\Delta^{-1,1}$. 
Thus $$\xi \wedge \eta = \sum_{\alpha <_t \beta} (\xi^\alpha \eta^\beta
- \xi^\beta \eta^\alpha) \cdot \ea \wedge \eb\;\;.$$
We deduce from the lemma~:
\begin{equation*} 
\begin{split}
\Xi(\xi \wedge \eta \wedge \ov{\xi \wedge \eta}) & = -
\sum_{\substack{ \alpha+ \beta = \gamma + \varepsilon \\\alpha <_t\beta \\ \gamma <_t \varepsilon \\ \alpha \not =
    \gamma  }} (\xi^\alpha \eta^\beta
- \xi^\beta \eta^\alpha) \cdot \ov{(\xi^\gamma \eta^\varepsilon-
  \xi^\varepsilon \eta^\gamma)} \cdot N_{\alpha, \beta} \cdot  N_{\gamma,
  \varepsilon} \\
& \qquad + \sum_{ \alpha <_t \beta} |\xi^\alpha \eta^\beta
- \xi^\beta \eta^\alpha|^2 \cdot (N_{\alpha, - \beta}^2 + <\alpha,
\beta>) 
\end{split}
\end{equation*}
and thus 
\begin{equation} \label{equa1}
\begin{split}
\Xi(\xi \wedge \eta \wedge \ov{\xi \wedge \eta}) 
& = - \sum_{\gamma <_t \varepsilon} \ov{(\xi^\gamma \eta^\varepsilon-
  \xi^\varepsilon \eta^\gamma)} \cdot N_{\gamma,
  \varepsilon}\cdot \left \{ \sum_{\substack{\alpha <_t\beta  \\ \alpha \not =
    \gamma \\ \alpha+ \beta = \gamma + \varepsilon}}  (\xi^\alpha \eta^\beta
- \xi^\beta \eta^\alpha) N_{\alpha, \beta}  \right \} \\
& \qquad+ \sum_{\alpha <_t \beta} |\xi^\alpha \eta^\beta
- \xi^\beta \eta^\alpha|^2 \cdot (N_{\alpha, - \beta}^2 + <\alpha, \beta>)
\;\;.
\end{split}
\end{equation}

Suppose now that $[\xi , \eta]= 0$. Thus $\sum_{\alpha <_t \beta}(\xi^\alpha \eta^\beta
- \xi^\beta \eta^\alpha) \cdot [e_{\alpha}, e_{\beta}] =0$ which is
\begin{equation*} 
\sum_{\substack{\alpha <_t \beta \\ \alpha+ \beta \; \text{fixed}}} (\xi^\alpha \eta^\beta
- \xi^\beta \eta^\alpha) \cdot N_{\alpha, \beta} =0 \;\;.
\end{equation*}
which implies for any $\gamma, \varepsilon \in \Delta$~':
\begin{equation} \label{equa2}
\sum_{\substack{\alpha <_t \beta \\ \alpha+ \beta = \gamma + \varepsilon
  \\ \alpha \not= \gamma }} (\xi^\alpha \eta^\beta
- \xi^\beta \eta^\alpha) \cdot N_{\alpha, \beta} = - (\xi^\gamma
\eta^\varepsilon - \xi^\beta \eta^\gamma) \cdot N_{\gamma, \varepsilon}\;\;.
\end{equation}

Replacing
$$\sum_{\substack{\alpha <_t \beta \\ \alpha+ \beta = \gamma + \varepsilon
  \\ \alpha \not= \gamma }} (\xi^\alpha \eta^\beta
- \xi^\beta \eta^\alpha) \cdot N_{\alpha, \beta} $$ 
by 
$$ - (\xi^\gamma
\eta^\varepsilon - \xi^\beta \eta^\gamma) \cdot N_{\gamma, \varepsilon}$$ in
equation~(\ref{equa1}) we obtain 
theorem~\ref{pont1}.

\end{proof}

\subsection{Proof of theorem~\ref{pont}}

By theorem~\ref{pont1} it is equivalent to showing that 
$$ N_{\alpha, - \beta}^2 + <\alpha, \beta> +
N_{\alpha, \beta}^2 \geq 0$$
for any roots $\alpha, \beta \in \Delta^{-1,1}$. 
Otherwise $<\alpha, \beta> <0$. But then $\alpha, \alpha+ \beta,
\cdots, \alpha + k \beta$ is the $\beta$-string through $\alpha$, with
$$k = -2 \frac{<\alpha, \beta>}{<\beta, \beta>}\;\;.$$
On the other hand by \cite[p.151]{helg}~:
$$N_{\alpha, \beta}^2 = \frac{k}{2} <\alpha, \alpha> = - <\alpha,
\beta> \;\;.$$
Thus
$$ N_{\alpha, - \beta}^2 + <\alpha, \beta> +
N_{\alpha, \beta}^2 = N_{\alpha, - \beta}^2 \geq 0\;\;.$$

\diam

\section{Proof of corollary~\ref{factor} and theorem~\ref{main}}

We first prove corollary~\ref{factor}.
Let $\gamma \in \Phi^c$ be the unique root whose
associated block $\FF_{\lieg^{-1,1}_{\CC, \gamma}} :=\bigoplus_{\alpha \in( \gamma +
  <\theta>) \cap \Delta_+} \lieg_ \alpha$ is contained in $\HHom(E^k,
E^{k-1})$. Let $\Psi^\gamma\subset \Psi$ be the connected component of
the vertex $\gamma$ in the subdiagram of the Dynkin diagram of $\Psi$ generated by $\Phi \cup
\{\gamma\}$. Let $\Delta^\gamma$ be the root system associated to
$\Psi^\gamma$ and $\Delta^\gamma_+ \subset \Delta_+$ the set of positive
roots in $\Delta^\gamma$. The root system $\delta^\gamma$ is of
Hermitian type, its unique positive non-compact root in $\Psi^\gamma$
is $\gamma$ and the set $(\gamma +<\theta>) \cap \Delta_+$ identifies with the set of positive
non-compact roots of this system. 

Suppose now $r_k=1$. This means that $\Psi^\gamma \setminus
\{\gamma\}$ is still connected and the roots in $(\gamma +<\theta>)
\cap \Delta_+$ are strictly ordered. Let $\alpha, \beta$ be two such
roots, without loss of generality we
can assume $\alpha \geq \beta$. Thus 
$$ <\alpha, \beta> = <\beta, \beta> +  <\alpha- \beta, \beta> \geq
<\beta, \beta> >0 \;\;,$$
as $\alpha -\beta$ is still a positive combination of positive roots
(and $\beta$ too).

This proves that for any $\xi, \eta \in \lieg^{-1,1}_{\CC, \gamma}$
satisfying $[\xi, \eta]=0$ one has 
$$\pi^* p_1(\FF) (\xi \wedge \eta \wedge \ov{\xi} \wedge \ov{\eta})
<0$$ 
except if $\xi^\alpha \eta^\beta = \xi^\beta \eta^\alpha$ for every
$\alpha, \beta \in (\gamma +<\theta>) \cap \Delta_+$, i.e. if $\xi$
and $\eta$ are proportional.

This immediately implies that if $\overline{f_\rho}^*\, p_1(\FF)$
vanishes identically on $X$ then $q_k \circ f_{\rho}: \ti{X} \lo
D_k$ has rank $1$, where $D_k$ denotes the open $G$-orbit in
$\check{D}_k$ image of $D$ under $q_k$. 

We now argue as in \cite[section 4.5]{KlKoMau}. We apply a theorem of Kaup
\cite[Theorem 6]{Gra} which states that there exists a normal complex
curve $\Sigma $ and holomorphic maps $h:\tilde M\lo \Sigma$,
$\psi:\Sigma\lo \check{D}_k$, with $h$ surjective, such that
$f_{\rho}=\psi\circ h$~: 
\begin{equation} \label{fac} 
\xymatrix{
\ti{X} \ar[r]^{f_{\rho}} \ar[d]_{h}& D
\ar[d]^{q_{k}} \\
\Sigma \ar[r]_\psi & D_k  \;\;.
}
\end{equation}
Which finishes the proof of corollary~\ref{factor}.

\sspace
To complete the proof of the main statement $b^2(\Ga) + b^4(\Ga) >0$
of theorem~\ref{main} one can assume we have
a factorization as in diagram~(\ref{fac}) (otherwise $b^4(\Ga) >0$). We continue to argue
as in \cite[section 4.5]{KlKoMau}. As the map
$q_k \circ f_\rho$ is superhorizontal and $D_k$ admits a function
$\phi$ whose Levi form $L(\phi)$ is positive on the superhorizontal
tangent space to $D_k$ (c.f. \cite[prop.4.2]{KlKoMau}) the curve
$\Sigma$ is not $\proj^1\CC$, in particular $\Sigma$ is aspherical.
By \cite[prop.4.2]{KlKoMau} the
restriction to $T_{sh}D_{k}$ of the $G$-equivariant curvature form of
the canonical line bundle $K_{D_{k}}$ is positive thus its pull-back
under $q_k \circ f_{\rho}$ is a positive $(1,1)$-form on $X$ defining
a non-trivial element in $H^2(X, \RR)$. As $q_k \circ f_\rho$
factorizes through $\Sigma$ which is aspherical, this element kills
$\pi_2(X)$ thus belongs to the subspace $H^2(\Ga, \RR) \hookrightarrow
H^2(X, \RR)$. Which finishes the proof of the main statement of
theorem~\ref{main}.

For the statement concerning $H^4_b(\Ga, \RR)$ it is enough to notice
that by \cite{gromov} the continuous cohomology of $G$ lifts to
bounded cohomology. The commutative diagram
$$ 
\xymatrix{
H^\bullet_{cb}(G, \RR) \ar@{->>}[r] \ar[d]_{\rho^*} &
H^\bullet_{ct}(G, \RR) \ar[d]^{\rho^{*}} \\
H^\bullet_b(\Ga, \RR)  \ar[r] & H^\bullet(\Ga, \RR)}\;\;$$ implies the result.

\diam

\section{Proof of theorem~\ref{linear2} and theorem~\ref{lattice}}

\subsection{Proof of theorem~\ref{linear2}}
Suppose $\Ga$ as in theorem~\ref{linear2}. Let $\rho : \Ga \hookrightarrow GL(n,
\CC)$ be an unbounded embedding of $\Ga$. As $b^2(\Ga) =0$ and $\Ga$
is an infinite K\"ahler group the group $\Ga$ is schematically rigid
\cite[cor.3]{kling}. Arguing as in \cite{kling} we can assume that $X$
is the fundamental group of a smooth projective complex variety $X$. 

By proposition~\ref{prop1} $\rho$ is the
monodromy of a \cvhs\- 
$$
\xymatrix{
\ti{X} \ar[r]^{f_{\rho}} \ar[dr]_{\overline{f_{\rho}}} &G/V
\ar[d]^{\pi} \\
& G/K\;\;,}$$
where $\G$ is a real form of Hodge type of the Zariski-closure of $\Ga$ in $\GL(n,
\CC)$. The group $G$ is non compact as $\Ga$ is unbounded in $GL(n,
\CC)$. By replacing if necessary this \cvhs\- by a convenient
exterior power as in corollary~\ref{factor} and noticing that the
corresponding monodromy still defines a linear unbounded embedding of
$\rho$ we obtain from corollary~\ref{factor} and its proof that the
map 
\begin{equation} \label{pull}
c^* : H^4(\Ga, \RR)  \lo H^4(X, \RR)
\end{equation}
deduced from the classyfing
map $c: X \lo B\Ga$ is non-zero~: it maps the non-zero class $\rho^*
p_1(\FF) \in H^4(\Ga, \RR)$ to the non-zero De Rham class
$[\overline{f_\rho}^*\, p_1(\FF)] \in H^{2,2}(X, \RR)$.

Using Lefschetz's hyperplane theorem we can without loss of generality
assume that $X$ is a (connected) smooth projective complex surface. In
this case equation~(\ref{pull}) is equivalent to saying that 
\begin{equation} \label{push}
c_* : H_4(M, \RR) \simeq \RR \lo H_4(\Ga, \RR) 
\end{equation}
is injective.

We proceed by contradiction to prove theorem~\ref{linear2}. Suppose $b^4(\Ga) \leq 1$. By
\cite[theor.1]{kling} we can assume that $b^4(\Ga)
=1$. By the previous remarks we obtain that
$$c_* : H_4(M, \RR) \simeq \RR \stackrel{\sim}{\lo} H_4(\Ga, \RR)
\simeq \RR$$ 
is an isomorphism.

We now follow closely the proof of \cite[theor.1]{kling}. As
in \cite[section 5.1.1]{kling} we get that the universal cover $\ti{X}$
has the homotopy type of a $2$-dimensional CW-complex. But then 
the homological exact sequence~(\cite[proposition 2, (a)]{kling}
$$H_4(M)_\QQ \lo H_4(\Ga)_\QQ \lo H_1(\Ga, \pi_2(M)_\QQ) \lo H_3(M)_\QQ $$
still leads to
\begin{equation} 
H_1(\Ga, \pi_2(M)_\QQ)=0\;\;.
\end{equation}
Following the end of the proof of \cite[theor.5]{kling} we still get a
contradiction to the duality property for $\Ga$.

\diam

\subsection{Proof of theorem~\ref{lattice}}
Let $G$ be a simple real Lie group of non-compact type.
Matsushima proved in \cite{mat} that there exists a constant $m(G)$
depending only of the Lie algebra $\lieg$ such that the restriction map
$$ \forall i \leq m(G), \quad H^i_{ct}(G, \RR) \lo H^i(\Ga, \RR)$$ is
an isomorphism for any cocompact lattice $\Ga$ of $G$.
One easily checks \cite{mat} that $m(G) \geq \frac{\rk_\RR G }{4} -1$.

For $G$ of Hodge type but not of Hermitian type we have 
$$\dim_\RR H^2(\lieg, K, \RR)=0 \quad \text{and} \quad \dim_\RR H^4(\lieg, K, \RR)=1\;\;.$$ 
Thus for $\rk_\RR G \geq 12$ one has $b^2(\Ga) =0$ and for $\rk_\RR G
\geq 20$ one has $b^4(\Ga) =1$. 

Moreover $\Ga$ is an $r$-dimensional duality group, with $r = \dim_\RR
(G/K) \geq \rk_\RR G$. 

Thus theorem~\ref{linear2} implies theorem~\ref{lattice}.

\diam


\lspace
\noi
{\small Bruno Klingler

\noi Institut de Math{\'e}matiques de Jussieu, Paris 75013, France

\noi
e-mail~: klingler@math.jussieu.fr}

\end{document}